\documentclass[12pt]{article}
\usepackage{amssymb}
\usepackage{latexsym,bm}
\usepackage{graphicx}
\usepackage{amsmath}
\usepackage{mathrsfs}

\newcommand{\bea}{\begin{eqnarray*}}
\newcommand{\eea}{\end{eqnarray*}}
\newcommand{\be}{\begin{equation}}
\newcommand{\ee}{\end{equation}}
\newcommand{\ben}{\begin{eqnarray*}}
\newcommand{\een}{\end{eqnarray*}}

\date{}
\begin{document}
\title{The maximum size of a nonhamiltonian graph with given order and connectivity\footnote{E-mail addresses:
{\tt zhan@math.ecnu.edu.cn}(X.Zhan),
{\tt mathdzhang@163.com}(L.Zhang).}}
\author{\hskip -10mm Xingzhi Zhan and Leilei Zhang\thanks{Corresponding author.}\\
{\hskip -10mm \small Department of Mathematics, East China Normal University, Shanghai 200241, China}}\maketitle
\begin{abstract}
 Motivated by work of Erd\H{o}s, Ota determined the maximum size $g(n,k)$ of a $k$-connected nonhamiltonian graph of order $n$ in 1995.
 But for some pairs $n,k,$ the maximum size is not attained by a graph of connectivity $k.$ For example, $g(15,3)=77$ is attained
 by a unique graph of connectivity $7,$ not $3.$ In this paper we obtain more precise information by determining the maximum size
 of a nonhamiltonian graph of order $n$ and connectivity $k,$ and determining the extremal graphs. Consequently we solve the corresponding problem
 for nontraceable graphs.
\end{abstract}



\section{Introduction}

We consider finite simple graphs, and use standard terminology and notations. The {\it order} of a graph is its number of vertices, and the
{\it size} its number of edges. For graphs we will use equality up to isomorphism, so $G_1=G_2$ means that $G_1$ and $G_2$ are isomorphic.
$\overline{G}$ denotes the complement of a graph $G.$ For two graphs $G$ and $H,$ $G\vee H$
denotes the {\it join} of $G$ and $H,$ which is obtained from the disjoint union $G+H$ by adding edges joining every vertex of $G$ to every vertex of $H.$
$K_n$ denotes the complete graph of order $n.$

One way to understand hamiltonian graphs is to investigate nonhamiltonian graphs. In 1961 Ore [8] determined the maximum size of a nonhamiltonian graph
with a given order and also determined the extremal graphs.

{\bf Lemma 1.} (Ore [8]) {\it The maximum size of a nonhamiltonian graph of order $n$ is $\binom{n-1}{2}+1$ and this size is attained by a graph $G$
if and only if $G=K_1\vee (K_{n-2}+K_1)$ or $G=K_2\vee \overline{K_3}.$}

Bondy [1] gave a new proof of Lemma 1. It is natural to ask the same question by putting constraints on the graphs. In 1962 Erd\H{o}s [6] determined
the maximum size of a nonhamiltonian graph of order $n$ and minimum degree at least $k,$ while in 1995 Ota [9] determined the maximum size
$g(n,k)$ of a $k$-connected nonhamiltonian graph of order $n.$ But for some pairs $n,k,$ the maximum size is not attained by a graph of connectivity $k.$
For example, $g(15,3)=77$ is attained  by a unique graph of connectivity $7,$ not $3.$

In this paper we obtain more precise information by determining the maximum size of a nonhamiltonian graph of order $n$ and connectivity $k,$ and
determining the extremal graphs, from which Ota's result can be deduced. Consequently we solve the corresponding problem for nontraceable graphs.

\section{Main results}

Denote by $V(G)$ and $E(G)$ the vertex set and edge set of a graph $G$ respectively. For $S\subseteq V(G),$ we denote by $G[S]$ the subgraph of $G$ induced by $S.$ ${\rm deg}(v)$ denotes the  degree of a vertex $v,$ and $\delta(G)$ denotes the minimum degree of a graph $G.$
Let $K_{s,\, t}$ denote the complete bipartite graph on $s$ and $t$ vertices.

We denote by $\kappa (G)$ and $\alpha(G)$ the connectivity and independence number of a graph $G$ respectively.

We will need the following lemmas.

{\bf Lemma 2.} (Chv\'{a}tal [4]) {\it Let $G$ be a graph with degree sequence $d_1\le d_2\le\cdots\le d_n$ where $n\ge 3.$
If there is no integer $k$ with $1\le k<n/2$ such that $d_k\le k$ and $d_{n-k}<n-k,$ then $G$ is hamiltonian.}

Lemma 2 can also be found in [3, p.488].

{\bf Lemma 3.} (Chv\'{a}tal-Erd\H{o}s [5]) {\it Let $G$ be a graph of order at least three. If $\kappa (G)\ge \alpha(G),$ then $G$ is hamiltonian.}

Lemma 3 can also be found in [3, p.488] and [10, p.292].

Given a graph $G$ and a positive integer $s$ with $s\le \alpha (G),$  denote
$$
\sigma_s(G)={\rm min}\left\{\left.\sum_{v\in T}{\rm deg} (v)\right| \, T\subseteq V(G)\,\,\, {\rm is\,\,\,an\,\,\,independent\,\,\, set{}}\,\,\, {\rm and}
\,\,\, |T|=s\right\}.
$$

The following result is a special case of Ota's theorem.

{\bf Lemma 4.} (Ota [9, Theorem 1]) {\it Let $G$ be a $k$-connected graph of order $n$ where $2\le k< \alpha(G).$ If for every integer $p$ with
$k\le p\le\alpha(G)-1$ we have $\sigma_{p+1}(G)\ge n+p^2-p,$ then $G$ is hamiltonian. }

A bipartite graph with partite sets $X$ and $Y$ is called {\it balanced} if $|X|=|Y|.$ For $n\ge 3,$ we denote by $K_{n,\, n-2}+4e$ the bipartite graph
obtained from $K_{n,\, n-2}$ by adding two vertices which are adjacent to two common vertices of degree $n-2.$

{\bf Lemma 5.} (Liu-Shiu-Xue [7]) {\it Given an integer $n\ge 4,$ let $\Omega(n)$ denote the set of all nonhamiltonian balanced bipartite graphs
of order $2n$ with minimum degree at least $2,$ and let $\Omega(3)$ denote the set of all nonhamiltonian balanced bipartite graphs of order $6.$
Then for any $n\ge 3,$ the maximum size of a graph in $\Omega(n)$ is $n^2-2n+4$ and this maximum size is uniquely attained by the graph
$K_{n,\, n-2}+4e.$ }

The case $n\ge 4$ of Lemma 5 is proved in [7, p.257] and the case $n=3$ can be verified easily. We will use this lemma with all
$n\ge 3$ cases.

{\bf Lemma 6.} (Bondy [2]) {\it Let $G$ be a graph of order $n$ with degree sequence $d_1\le d_2\le\cdots\le d_n$ and let $k$ be an integer with $0\le k\le n-2.$
If for each integer $j$ with $1\le j\le n-1-d_{n-k}$ we have $d_j\ge j+k,$ then $G$ is $(k+1)$-connected. }

{\bf Lemma 7.} {\it Let $G=K_s\vee\overline{K_t}$ or $G=K_s\vee(K_2+\overline{K_t})$ where $t\ge 2$ in both cases, and let $F\subseteq E(G)$ with $|F|=f\le s.$
Then $\kappa (G-F)\ge s-f,$ with equality if and only if all the edges in $F$ are incident to one common vertex in $\overline{K_t}.$ }

{\bf Proof.} We prove the case when $G=K_s\vee\overline{K_t}.$ The case when $G=K_s\vee(K_2+\overline{K_t})$ can be proved similarly.

It is easy to see that $\kappa(G)=s.$ Since deleting one edge reduces the connectivity by at most one [10, p.169], we have $\kappa(G-F)\ge s-f.$

Next we use induction on $f$ to prove the equality condition. First consider the case $f=1.$ Let $e\in E(G).$ It is easy to check that $\kappa(G-e)=s-1$
if and only if $e$ has one endpoint in $K_s$ and the other endpoint in $\overline{K_t}.$ Now let $F\subseteq E(G)$ with $|F|=f\ge 2$ and suppose that
for any $A\subseteq E(G)$ with $|A|=f-1,$ $\kappa(G-A)=s-(f-1)$ if and only if all the edges in $A$ are incident to one common vertex in $\overline{K_t}.$

If all edges in $F$ are incident to one common vertex in $\overline{K_t},$ it is easy to verify that $\kappa(G-F)=s-f.$ Conversely, suppose
$\kappa(G-F)=s-f.$ Let $F=\{e_1,e_2,\ldots,e_f\}$ and denote $F^{\prime}=F\setminus\{e_f\}.$ Then $\kappa(G-F^{\prime})=s-f+1.$ By the induction hypothesis,
the edges $e_1,\ldots,e_{f-1}$  are incident to one common vertex $w$ in $\overline{K_t}.$ The degree sequence of $G-F^{\prime}$ is
$$
s-f+1,\, \underbrace{s,\ldots,s}_{t-1},\, \underbrace{n-2,\ldots,n-2}_{f-1},\, \underbrace{n-1,\ldots,n-1}_{s-f+1}
$$
where $n=s+t$ and $s-f+1={\rm deg}(w).$ We assert that $e_f$ is incident to $w$ and consequently all the edges in $F$ are incident to one common vertex in $\overline{K_t}.$ Let the degree sequence of $G-F$ be $d_1\le\cdots\le d_n.$ By the above degree sequence of $G-F^{\prime}$ we deduce that
$d_{n-s+f}\ge n-2.$ Thus $n-1-d_{n-s+f}\le 1.$ If $e_f$ is not incident to $w,$ then we would have $d_1=s-f+1=1+(s-f).$ By Lemma 6, $G-F$ is
$(s-f+1)$-connected, contradicting the assumption $\kappa(G-F)=s-f.$ This proves that  $e_f$ is incident to $w.$ \hfill $\Box$

{\bf Notation 1.} $e(G)$ denotes the size of a graph $G.$

{\bf Notation 2.} For positive integers $n$ and $k$ with $n$ odd and $n\ge 2k+1,$ $G_1(n,k)$ denotes the graph obtained from
$K_{(n-1)/2}\vee\overline{K_{(n+1)/2}}$ by deleting $(n-1)/2-k$ edges that are incident to one common vertex in $\overline{K_{(n+1)/2}};$
for positive integers $n$ and $k$ with $n$ even and $n\ge 2k+2,$ $G_2(n,k)$ denotes the graph obtained from
$K_{(n-2)/2}\vee(K_2+\overline{K_{(n-2)/2}})$ by deleting $(n-2)/2-k$ edges that are incident to one common vertex in $\overline{K_{(n-2)/2}}.$

Note that by Dirac's theorem [3, p.485], for the existence of a nonhamiltonian graph of order $n$ and connectivity $k$ we necessarily
have $n\ge 2k+1.$ Now we are ready to state and prove the main result.

{\bf Theorem 8.} {\it Let $f(n,k)$ denote the maximum size of a nonhamiltonian graph of order $n$  and connectivity $k.$ Then
$$
f(n,k)=\begin{cases} \binom{n-k}{2}+k^2\quad {\rm if}\,\,\,n\,\,\,{\rm is}\,\,\,{\rm odd}\,\,\,{\rm and}\,\,\,n\ge 6k-5\,\,\,{\rm or}\,\,\,
n\,\,\,{\rm is}\,\,\,{\rm even}\,\,\,{\rm and}\,\,\,n\ge 6k-8,\\
\frac{3n^2-8n+5}{8}+k\quad {\rm if}\,\,\,n\,\,\,{\rm is}\,\,\,{\rm odd}\,\,\,{\rm and}\,\,\,2k+1\le n\le 6k-7,\\
\frac{3n^2-10n+16}{8}+k\quad {\rm if}\,\,\,n\,\,\,{\rm is}\,\,\,{\rm even}\,\,\,{\rm and}\,\,\,2k+2\le n\le 6k-10.
\end{cases}
$$
If $n=6k-5,$ then $f(n,k)$ is attained by a graph $G$ if and only if $G=K_k\vee(K_{n-2k}+\overline{K_k})$ or $G=G_1(n,k).$
If $n=6k-8,$ then $f(n,k)$ is attained by a graph $G$ if and only if $G=K_k\vee(K_{n-2k}+\overline{K_k})$ or $G=G_2(n,k).$
If $n$ is odd and $n\ge 6k-3$ or $n$ is even and $n\ge 6k-6,$ then $f(n,k)$ is attained by a graph $G$ if and only if $G=K_k\vee(K_{n-2k}+\overline{K_k}).$
If $n$ is odd and $2k+1\le n\le 6k-7,$ then $f(n,k)$ is attained by a graph $H$ if and only if $H=G_1(n,k).$
If $n$ is even and $2k+2\le n\le 6k-10,$ then $f(n,k)$ is attained by a graph $Z$ if and only if $Z=G_2(n,k).$ }

{\bf Proof.} The case $k=1$ of Theorem 8 follows from Lemma 1. Note that the extremal graph $K_2\vee \overline{K_3}$ of order $5$ in Lemma 1 has
connectivity $2$ and hence it should be excluded.

Next suppose $k\ge 2.$ It is easy to verify that the extremal graphs stated in Theorem 8 are nonhamiltonian graphs of order $n$ and connectivity $k$
with size $f(n,k).$ They are nonhamiltonian since any hamiltonian graph must be tough [3, pp.472-473]. Thus it remains to show that $f(n,k)$
is an upper bound on the size and it can only be attained by these extremal graphs.

Let $Q$ be a nonhamiltonian graph of order $n$ and connectivity $k$ with degree sequence $d_1\le d_2\le\cdots\le d_n.$ By Lemma 3, $k<\alpha(Q)$
and by Lemma 4, there exists an integer $p$ with $k\le p\le \alpha(Q)-1$ such that $\sigma_{p+1}(Q)\le n+p^2-p-1.$ Let $S$ be an independent set of $Q$ with cardinality $p+1$
whose degree sum is $\sigma_{p+1}(Q).$ Then $e(Q[V(Q)\setminus S])\le \binom{n-p-1}{2}.$ We distinguish four cases.

Case 1. $n$ is odd and $n\ge 6k-5.$

Subcase 1.1. $p\le (n-3)/2.$

The conditions  $p\le (n-3)/2$ and $n\ge 6k-5$ imply $3p+3k+1<2n.$ This, together with the condition $p\ge k,$ yields
$(p-k)(3p+3k+1-2n)\le 0.$ It follows that
$$
e(Q)\le n+p^2-p-1+\binom{n-p-1}{2}\le \binom{n-k}{2}+k^2 \eqno (1)
$$
and equality holds in the second inequality in (1) if and only if $p=k.$

Now suppose that $Q$ has size $\binom{n-k}{2}+k^2.$ Then $p=k,$ $S$ has cardinality $k+1$ and degree sum $n+k^2-k-1,$ and
$V(Q)\setminus S$ is a clique. Since $k+1<(n+1)/2,$ we have $d_{(n+1)/2}\ge n-k-2.$ By Lemma 2, there exists $i$ with $i<n/2$ such that $d_i\le i$ and $d_{n-i}\le n-i-1.$
Since $n$ is odd, the condition $i<n/2$ means $i\le (n-1)/2.$ We have
$$
e(Q)=\binom{n-k}{2}+k^2\le [i^2+(n-2i)(n-i-1)+i(n-1)]/2, \eqno (2)
$$
where the inequality is equivalent to $(i-k)(2n-3i-3k-1)\le 0.$ Since $i\ge d_i\ge \delta(Q)\ge k,$ we obtain $i=k$ or $n\le (3i+3k+1)/2.$

If $i=k,$ equality holds in (2) and hence the degree sequence of $Q$ is
$$
\underbrace{k,\ldots,k}_{k},\, \underbrace{n-k-1,\ldots,n-k-1}_{n-2k},\, \underbrace{n-1,\ldots,n-1}_{k},
$$
implying that $Q=K_k\vee(K_{n-2k}+\overline{K_k}).$

Now suppose $i\neq k.$ Then we have $n\le (3i+3k+1)/2.$ If $i\le (n-3)/2,$ then $n\le 6k-7,$ contradicting our assumption $n\ge 6k-5.$ Thus $i=(n-1)/2.$
We have $n-k-2\le d_{(n+1)/2}\le (n-1)/2.$ Hence $6k-5\le n\le 2k+3,$ which, together with the condition $k\ge 2,$ yields  $k=2$ and $n=7.$
It is easy to check that there are exactly four graphs of order $7$ and size $14$
with $d_4=3,$ among which $G_1(7,2)$ is the only graph that is nonhamiltonian with connectivity $2.$ Hence $Q=G_1(7,2).$

Subcase 1.2. $p\ge (n-1)/2.$

Clearly $Q$ is a spanning subgraph of $R=K_{n-p-1}\vee\overline{K_{p+1}}.$ If $p\ge (n+1)/2,$ then
$$
e(Q)\le \binom{n-p-1}{2}+(n-p-1)(p+1)<\binom{n-k}{2}+k^2
$$
where the second inequality follows from the condition $n\ge 6k-5.$ If $p=(n-1)/2,$ we have $\kappa(R)=n-p-1>k.$ Let $F\subseteq E(R)$ such that
$Q=R-F.$ Since $\kappa(Q)=k,$ by Lemma 7 we have $|F|\ge n-p-1-k.$ Thus
$$
e(Q)\le\binom{n-p-1}{2}+(n-p-1)(p+1)-(n-p-1-k)\le \binom{n-k}{2}+k^2 \eqno (3)
$$
and equality holds in the second inequality of (3) if and only if $n=6k-5.$

Suppose $e(Q)=\binom{n-k}{2}+k^2.$ Then $n=6k-5$ and $|F|=n-p-1-k.$ By Lemma 7, all the edges in $F$ are incident to one common vertex in
$\overline{K_{p+1.}}$ Since $p=(n-1)/2,$ we have $n-p-1=(n-1)/2,$ $p+1=(n+1)/2$ and $|F|=(n-1)/2-k.$ It follows that $Q=G_1(n,k).$

Case 2. $n$ is odd and $2k+1\le n\le 6k-7.$

Subcase 2.1. $k\le p< (n-1)/2.$

We have
$$
e(Q)\le n+p^2-p-1+\binom{n-p-1}{2}<  \frac{3n^2-8n+5}{8}+k.
$$

Subcase 2.2. $p=(n-1)/2.$

In this case $Q$ is a spanning subgraph of $K_p\vee\overline{K_{p+1}}.$ Since $\kappa(Q)=k,$ by Lemma 7 we obtain
$$
e(Q)\le \binom{p}{2}+p(p+1)-(p-k)=\frac{3n^2-8n+5}{8}+k
$$
and equality holds if and only if $Q=G_1(n,k).$

Subcase 2.3. $p>(n-1)/2.$

In this case $Q$ is a spanning subgraph of $K_{n-p-1}\vee\overline{K_{p+1}}.$ Then
$$
e(Q)\le\binom{n-p-1}{2}+(n-p-1)(p+1)<\frac{3n^2-8n+5}{8}+k,
$$
where we have used the condition $n\le 2p-1.$

Case 3. $n$ is even and $n\ge 6k-8.$

Subcase 3.1. $p<(n-2)/2.$

The assumptions imply $3p+3k+1-2n< 0.$ We have
$$
e(Q)\le n+p^2-p-1+\binom{n-p-1}{2}\le \binom{n-k}{2}+k^2  \eqno (4)
$$
where the second inequality is equivalent to
$$
(p-k)(3p+3k+1-2n)\le 0.
$$
Thus equality holds in the second inequality in (4) if and only if $p=k.$

Suppose $e(Q)=\binom{n-k}{2}+k^2.$ Then $p=k$ and $Q$ has a clique of cardinality $n-p-1$ and an independent set of cardinality $p+1$ whose degree sum
equals $n+p^2-p-1.$ Also $d_{(n+2)/2}\ge n-k-2.$ By Lemma 2, there exists $i<n/2$ such that $d_i\le i$ and $d_{n-i}\le n-i-1.$ We have
$$
e(Q)=\binom{n-k}{2}+k^2\le [i^2+(n-2i)(n-i-1)+i(n-1)]/2,
$$
where the inequality is equivalent to
$$
(i-k)(2n-3i-3k-1)\le 0. \eqno (5)
$$

Note that $i\ge k$ since $i\ge d_i\ge \delta(Q)\ge k.$
If $i=k,$ then the degree sequence of $Q$ is
$$
\underbrace{k,\ldots,k}_{k},\, \underbrace{n-k-1,\ldots,n-k-1}_{n-2k},\, \underbrace{n-1,\ldots,n-1}_{k},
$$
implying $Q=K_k\vee (K_{n-2k}+\overline{K_k}).$

Next suppose $i\neq k.$ Then the inequality (5) implies $2n-3i-3k-1\le 0.$ If $i\le (n-4)/2,$ we deduce that $n\le 6k-10,$ a contradiction.
Hence $i=(n-2)/2.$ Now the conditions $d_{n-i}\le n-i-1$ and $d_{(n+2)/2}\ge n-k-2$ yield $n-k-2\le d_{(n+2)/2}\le n/2.$
Thus $6k-8\le n\le 2k+4.$ It follows that $k=2$ and $4\le n\le 8$ or $k=3$ and $10\le n\le 10.$ The possibility $n=4$ contradicts $k<n/2$
and $n=6$ contradicts $i\neq k.$ Only the two pairs $(k,n)=(2,\, 8),\,(3,\, 10)$ can occur.

If $k=2$ and $n=8,$ by the conditions $n=8,$ $e=19,$ $d_5=4$ and being nonhamiltonian we deduce that $Q=K_3\vee(K_2+\overline{K_3});$ if $k=3$ and $n=10,$
the conditions $n=10,$ $e(Q)=30,$ $k=3,$ $d_6=5$ and being nonhamiltonian force $Q=G_2(10,3).$

Subcase 3.2. $p=(n-2)/2.$

Clearly $\alpha(Q)\ge p+1.$ We further distinguish two cases.

If $\alpha(Q)\ge p+2,$ then $Q$ is a spanning subgraph of $K_p\vee \overline{K_{p+2}}.$ By Lemma 7 we have
$$
e(Q)\le \binom{p}{2}+p(p+2)-(p-k)< \binom{n-k}{2}+k^2. \eqno (6)
$$
The second inequality in (6) is equivalent to $p^2+(5-4k)p+3k^2-5k+2 > 0$ which is guaranteed by $p=(n-2)/2\ge 3k-5.$

If $\alpha(Q)=p+1,$ then $Q$ is a spanning subgraph of $K_{p+1}\vee \overline{K_{p+1}}.$ Let $Q^{\prime}$ denote the graph obtained from $Q$ by deleting all the edges in $K_{p+1}.$ Then $Q^{\prime}$ is a nonhamiltonian balanced bipartite graph. There are two cases.

(a) Suppose $n\ge 8$ and $\delta (Q^{\prime})\ge 2$ or $n=6.$ By Lemma 5, $e(Q^{\prime})\le (p+1)^2-2(p+1)+4=p^2 +3.$ Hence
$$
e(Q)\le \binom{p+1}{2}+p^2+3\le \binom{n-k}{2}+k^2. \eqno(7)
$$
The second inequality in (7) is equivalent to
$$
p^2+(5-4k)p+3k^2-3k-4\ge 0 \eqno(8)
$$
which is implied by the condition $p\ge 3k-5.$ Equality holds in (8) if and only if $k=2$ and $p=2,$ i.e., $Q^{\prime}$ is the extremal graph of order $6$ defined
in Lemma 5. Hence $Q$ has size $\binom{n-k}{2}+k^2$ if and only if $Q=K_2\vee(K_2+\overline{K_2}).$

(b) Now suppose $n\ge 8$ and $\delta (Q^{\prime})\le 1.$ Let $x\in V(Q^{\prime})$ with ${\rm deg}_{Q^{\prime}}(x)=\delta (Q^{\prime}).$
Starting with the structure $K_{p+1}\vee \overline{K_{p+1}},$ we deduce that $x$ lies in $K_{p+1},$ since $\delta (Q)\ge 2.$
In this case $Q$ is a spanning subgraph of $K_p\vee (K_2+\overline{K_p}).$ By Lemma 7
and using the fact that $p\ge 3k-5$ we have
$$
e(Q)\le \binom{p+2}{2}+p^2-(p-k)\le \binom{n-k}{2}+k^2. \eqno(9)
$$
Equality in the second inequality in (9) holds if and only if $p=3k-5;$ i.e, $n=6k-8.$ Thus by (9) and Lemma 7, $Q$ has size $\binom{n-k}{2}+k^2$ if and only if
$Q=G_2(n,k)$ with $n=6k-8.$

Subcase 3.3. $p>(n-2)/2.$

$Q$ is a spanning subgraph of $K_{n-p-1}\vee \overline{K_{p+1}}.$ If $p=n/2,$ we have $n-p-1\ge k.$ By Lemma 7
$$
e(Q)\le \binom{n-p-1}{2}+(n-p-1)(p+1)-(n-p-1-k)< \binom{n-k}{2}+k^2. \eqno(10)
$$
The second inequality in (10) is equivalent to $p^2+(3-4k)p+3k^2-k-2>0,$ which is implied by $p=n/2\ge 3k-4.$

If $p\ge (n+2)/2,$ we have
$$
e(Q)\le \binom{n-p-1}{2}+(n-p-1)(p+1)< \binom{n-k}{2}+k^2. \eqno(11)
$$
The second inequality in (11) is equivalent to $p^2+p+3k^2+k-2nk>0.$ To prove this inequality it suffices to show
$p^2+(1-4k)p+3k^2+5k>0,$ which is implied by $p\ge (n+2)/2\ge 3k-3.$

Case 4. $n$ is even and $2k+2\le n\le 6k-10.$

Denote $m=(n-2)/2.$ Then $3\le k\le m.$ We distinguish three subcases.

Subcase 4.1. $p<m.$

We have
$$
e(Q)\le \binom{n-p-1}{2}+n+p^2-p-1< \frac{3n^2-10n+16}{8}+k. \eqno(12)
$$
The second inequality in (12) is equivalent to $3p^2-(4m+3)p+m^2+5m-2k<0,$ which is implied by the conditions $k\le p<m\le 3k-6.$

Subcase 4.2. $p=m.$

If $\alpha(Q)\ge p+2,$ then $Q$ is a spanning subgraph of $K_p\vee \overline{K_{p+2}}.$ Recall that $p\ge k.$ By Lemma 7 we have
\begin{align*}
e(Q)\le \binom{p}{2}+p(p+2)-(p-k)=(3p^2+p+2k)/2&<(3p^2+p+2k)/2+1\\
                                               &=\frac{3n^2-10n+16}{8}+k.
\end{align*}

Now suppose $\alpha(Q)=p+1.$ Then $Q$ is a spanning subgraph of $K_{p+1}\vee \overline{K_{p+1}}.$
Define the graph $Q^{\prime}$ as in Subcase 3.2 above.

If $\delta(Q^{\prime})\ge 2,$ by Lemma 5 we have $e(Q^{\prime})\le p^2+3.$ Hence
$$
e(Q)\le \binom{p+1}{2}+p^2+3<\frac{3n^2-10n+16}{8}+k.
$$

If $\delta(Q^{\prime})\le 1,$ then $Q$ is a spanning subgraph of $K_p\vee (K_2+\overline{K_p}).$ By Lemma 7 we obtain
$$
e(Q)\le \binom{p+2}{2}+p^2-(p-k)=\frac{3n^2-10n+16}{8}+k
$$
and equality holds if and only if $Q=G_2(n,k).$

Subcase 4.3. $p>m.$

$Q$ is a spanning subgraph of $K_{n-p-1}\vee \overline{K_{p+1}}.$ We further distinguish two cases.

(a) $p=m+1.$ Now the conditions $2k+2\le n$ and $p=m+1=n/2$ imply $n-p-1=(n-2)/2\ge k.$ By Lemma 7, we have
$$
e(Q)\le \binom{n-p-1}{2}+(n-p-1)(p+1)-(n-p-1-k)< \frac{3n^2-10n+16}{8}+k
$$
where the second inequality is equivalent to $4p^2-n^2+2n-4p+8>0$ which holds, since $n=2p.$

(b) $p\ge m+2.$ In this case the following rough estimate suffices:
$$
e(Q)\le \binom{n-p-1}{2}+(n-p-1)(p+1)<\frac{3n^2-10n+16}{8}+k. \eqno(13)
$$
The second inequality in (13) is equivalent to $4p^2+4p-n^2-6n+8k+16>0,$ which holds, since $n\le 2p-2.$
This completes the proof. \hfill $\Box$

The following corollary follows from Theorem 8 immediately.

{\bf Corollary 9.} {\it Let $f(n,k)$ be defined as in Theorem 8. If $G$ is a graph of order $n$ and connectivity $k$ with size larger than
$f(n,k),$ then $G$ is hamiltonian.}

Next we use Theorem 8 to deduce Ota's result.

{\bf Corollary 10.} (Ota [9, p.209]) {\it The maximum size of a  $k$-connected nonhamiltonian graph of order $n$ is
$$
{\rm max}\left\{\binom{n-k}{2}+k^2,\,\, \binom{\lfloor(n+2)/2\rfloor}{2}+\left\lfloor\frac{n-1}{2}\right\rfloor ^2\right\}. \eqno(14)
$$
}

{\bf Proof.} Denote the number in (14) by $M.$ Let $f(n,c)$ be defined as in Theorem 8.
It is easy to verify that
$$
{\rm max}\{f(n,c)|\,\, k\le c< n/2\}=M
$$
and the result follows. \hfill $\Box$

A graph is {\it traceable} if it contains a Hamilton path; otherwise it is {\it nontraceable.} Next we turn to nontraceable graphs.

The following trick is well-known (e.g. [4, p.166] or [5, p.112]).

{\bf Lemma 11.} {\it Let $G$ be a graph and denote $H=G\vee K_1.$ Then $G$ is traceable if and only if $H$ is hamiltonian, and
$\kappa(G)=k$ if and only if $\kappa (H)=k+1.$ }

{\bf Notation 3.} For positive integers $n$ and $k$ with $n$ odd and $n\ge 2k+3,$ $H_1(n,k)$ denotes the graph obtained from
$K_{(n-3)/2}\vee (K_2+\overline{K_{(n-1)/2}})$ by deleting $(n-3)/2-k$ edges that are incident to one common vertex in $\overline{K_{(n-1)/2}};$
for positive integers $n$ and $k$ with $n$ even and $n\ge 2k+2,$ $H_2(n,k)$ denotes the graph obtained from
$K_{(n-2)/2}\vee\overline{K_{(n+2)/2}}$ by deleting $(n-2)/2-k$ edges that are incident to one common vertex in $\overline{K_{(n+2)/2}}.$

By Dirac's theorem [3, p.485] and Lemma 11, for the existence of a nontraceable graph of order $n$ and connectivity $k$ we must
have $n\ge 2k+2.$ The following corollary follows from Theorem 8 and Lemma 11 immediately.

{\bf Corollary 12.} {\it Let $\varphi(n,k)$ denote the maximum size of a nontraceable graph of order $n$  and connectivity $k.$   Then
$$
\varphi(n,k)=\begin{cases} \binom{n-k-1}{2}+k(k+1)\quad {\rm if}\,\,\,n\,\,\,{\rm is}\,\,\,{\rm odd}\,\,\,{\rm and}\,\,\,n\ge 6k-3\,\,\,{\rm or}\,\,\,
n\,\,\,{\rm is}\,\,\,{\rm even}\,\,\,{\rm and}\,\,\,n\ge 6k,\\
\frac{3n^2-12n+17}{8}+k\quad {\rm if}\,\,\,n\,\,\,{\rm is}\,\,\,{\rm odd}\,\,\,{\rm and}\,\,\,2k+3\le n\le 6k-5,\\
\frac{3n^2-10n+8}{8}+k\quad {\rm if}\,\,\,n\,\,\,{\rm is}\,\,\,{\rm even}\,\,\,{\rm and}\,\,\,2k+2\le n\le 6k-2.
\end{cases}
$$
If $n=6k-3,$ then $\varphi(n,k)$ is attained by a graph $G$ if and only if $G=K_k\vee(K_{n-2k-1}+\overline{K_{k+1}})$ or $G=H_1(n,k).$
If $n=6k,$ then $\varphi(n,k)$ is attained by a graph $G$ if and only if $G=K_k\vee(K_{n-2k-1}+\overline{K_{k+1}})$ or $G=H_2(n,k).$
If $n$ is odd and $n\ge 6k-1$ or $n$ is even and $n\ge 6k+2,$ then $\varphi(n,k)$ is attained by a graph $G$ if and only if $G=K_k\vee(K_{n-2k-1}+\overline{K_{k+1}}).$
If $n$ is odd and $2k+3\le n\le 6k-5,$ then $\varphi(n,k)$ is attained by a graph $G$ if and only if $G=H_1(n,k).$
If $n$ is even and $2k+2\le n\le 6k-2,$ then $\varphi(n,k)$ is attained by a graph $G$ if and only if $G=H_2(n,k).$ }

\vskip 5mm
{\bf Acknowledgement.} This research  was supported by the NSFC grants 11671148 and 11771148 and Science and Technology Commission of Shanghai Municipality (STCSM) grant 18dz2271000.

\end{document}